\documentclass[12pt,letterpaper,fleqn]{article}

\title{A robust optimization model for green supplier selection and order allocation in a closed-loop supply chain considering cap-and-trade mechanism}

\author{
  Hossein Mirzaee$^{1,2}$\\
  \texttt{hossein.mirzaee@usask.ca}
  \and
  Hamed Samarghandi$^{2}$\\
  \texttt{samarghandi@edwards.usask.ca}
  \and
  Keith Willoughby$^{2}$\\
  \texttt{willoughby@edwards.usask.ca}\blfootnote{1 Corresponding author}\blfootnote{2 Edwards School of Business, University of Saskatchewan, Saskatoon, SK, Canada}
}

\usepackage{lipsum}

\newcommand\blfootnote[1]{%
  \begingroup
  \renewcommand\thefootnote{}\footnote{#1}%
  \addtocounter{footnote}{-1}%
  \endgroup
}

\usepackage[natbibapa]{apacite}
\usepackage{cancel}
\usepackage[toc,page,title]{appendix}
\usepackage{adjustbox}
\usepackage[dvipsnames]{xcolor}
\usepackage{comment}
\usepackage{hyperref}
\hypersetup{colorlinks,linkcolor={blue},citecolor={blue},urlcolor={blue}}  
\usepackage{amsmath}
\usepackage[left=1in,right=1in,top=1in,bottom=1in]{geometry}
\usepackage{xcolor}
\usepackage{booktabs}
\usepackage{float}
\usepackage{tikz}
\usepackage{pgfplotstable}
\usepackage[normalem]{ulem}
\usepackage{pgfplots}
\usepackage{url}
\usepackage{graphicx}
\usepackage[justification=centering]{caption}
\usetikzlibrary{positioning,calc}       
 \tikzstyle{block} = [draw,rectangle,thick,minimum height=2em,minimum width=2em,fill=green!10!white]
 \def\checkmark{\tikz\fill[scale=0.4](0,.35) -- (.25,0) -- (1,.7) -- (.25,.15) -- cycle;} 
\usepackage[utf8x]{inputenc}
\usepackage{graphicx}
\usepackage[english]{babel}
\date{}

\begin{document}

\maketitle
\begin{abstract}
\noindent Due to increasing air pollution, which is a consequence of the environmental effects of production in various industries, green supply chain management (GSCM) has attracted the attention of both scholars and practitioners. Green supplier selection is one of the most important problems in GSCM, which satisfies a firm's environmental goals as well as its economic targets. In this paper, for the first time, a green supplier selection problem considering both green and non-green evaluation criteria in a closed-loop supply chain is studied, and a cap-and-trade mechanism as a way of controlling the air pollution caused by manufacturers is proposed. To solve the described problem, a multi-objective robust optimization (RO) model, as an effective approach to handle uncertainty, is proposed. A numerical example using randomly generated data, accompanied by the analysis based on the proposed approach is elaborated to validate the presented model. The results prove that the developed model for green supplier selection is able to effectively enhance the decision-making process of the experts. By illustrating the trade-off in robustness between the model and proposed solutions, as well as the effect of the deviation penalty on the closeness of results to the achieved solution, we show how firms can make optimal decisions when assigning the parameters. Furthermore, analyses show that allowance amount (cap) and allowance prices in the cap-and-trade system impact the firms' costs and the amount of carbon released. Finally, we show that the cap-and-trade mechanism results in a better solution in terms of the total utility of the supply chain compared to the penalty-based system.

\vspace{1cm}
\end{abstract}

\setlength{\columnsep}{0.6cm}

\section{Introduction}\label{Introduction}
Due to the increasing environmental issues that have arisen recently, GSCM has become an important topic for both practitioners and researchers \citep{sang2016interval}. According to \cite{3Srivastava2007}, GSCM brings environmental considerations to supply chain management; this includes product design, raw material and parts sourcing, production processes, and transportation of finished products to customers. One of the most critical topics in GSCM is supplier selection, since about 70\% of the cost of the final product arise from component parts and raw material \citep{ghodsypour2001total}. Thus, procuring the raw materials from proper suppliers significantly impacts the characteristics of the final product. Accordingly, companies bestow privilege to the partners that score the highest on  environmental criteria \citep{3rao2005}.

Supplier selection is the process of evaluating various suppliers to select the best one(s) according to the high-priority criteria for the manufacturers \citep{ding2015decision}. In case of multiple sourcing, where more than one partner can be selected for outsourcing, a solution must encompass the selected suppliers and the size of the order to be placed with each of them. To evaluate the suppliers, the manufacturers should first specify their appraisal criteria. According to the literature of GSCM, there are two general types of measures in supplier selection: non-green and green.

Non-green criteria ascertain a company's competitiveness compared to other firms in the market. Among the non-green measures used in the selection process in the literature, the most common ones are: cost, quality, and lead time \citep{35}. In contrast, green criteria contemplate those of a company's activities with an adverse impact on the quality of water, air, and soil. The most utilized measures in evaluating the green criteria are the toxicity level of the materials used in products, recyclability, green production, and environmental management systems \citep{8}, as well as pollution production (i.e., the amount of carbon emission) \citep{10}.

Considering green criteria mentioned above by manufacturers is caused by some incentives that government makes for them to reduce the pollution level. One of these incentives can be restricting the emission level in a way that going beyond that level is costly for polluters. Carbon emission is usually measured by production level as well as equipment used in the supply chain \citep{10}. Companies usually produce products to maximize their profit; they do not supervise the amount of carbon emission during the production process. Recently, a mechanism called cap-and-trade has been proposed by various governments to keep the level of $CO_{2}$ emissions at a reasonable level. In this mechanism, a cap is defined as the maximum carbon emission allowance for manufacturers, and is treated as the emission quota when producing the products and responding to the customer demands. Furthermore, companies can sell or buy carbon emission allowances from each other in a trading market \citep{xu2022region}. In this situation, companies alter their objective to profit maximization not only by manufacturing the proper amount of products, but also by buying or selling carbon emission allowances. Generally, in the cap-and-trade mechanism, the assigned cap is reduced each year by the government to reach an acceptable level of carbon pollution.

In a supply chain, the amount of production, and consequently, the decisions related to raw material outsourcing are highly impacted by market demand and procurement cost. Typically, these factors are not under the companies' control. Although firms can impact these parameters by marketing and product pricing, they cannot predict their exact values because of unpredictable factors such as rivals' strategies. This uncertainty prevents the manufacturer from forecasting precise values for demand and procurement cost. Also, in a closed-loop supply chain, the quantity of returns consisting of used or rejected products is not easy to predict and should be considered uncertain \citep{1}. Consequently, a significant challenge for decision-makers is handling the uncertainty of the parameters such that their proposed solutions represent the real-life scenarios. In the context of suply chain, there are three common ways to control the uncertainty of the input data: stochastic optimization, fuzzy set theory, and robust optimization (RO) \citep{tordecilla2021simulation}.

The majority of studies available in the literature of green supplier selection has considered the uncertain parameters as stochastic or fuzzy. One of the challenges that the decision-makers face with stochastic programming is estimating the correct probability distribution of the uncertain parameters \citep{21}. However, in many cases, historical data related to the uncertain parameters is not available. Hence, the estimated probability distributions might not be reliable \citep{vahdani2012}. Furthermore, despite the high probability of correct prediction of the uncertain parameters, in very rare cases, it may be wrong and cause solution infeasibility. Even though this issue may occur with a very small probability, it can lead to a serious rise in the costs \citep{1}. What about fuzzy set theory? When using fuzzy set theory, exhaustive knowledge and comprehension about the parameters is required to generate an accurate membership function. Undoubtedly, obtaining such comprehensive consciousness about the market and its dynamics poses an enormous challenge for decision-makers, and is considered one of the obstacles of effectively employing the fuzzy set theory for real-world problem optimization. \citep{memon2015}. 

Robust optimization is a functional technique to handle uncertainties, and is mostly free from the flaws mentioned above \citep{jabbarzadeh2019}. Both fuzzy optimization and stochastic optimization rely on probability distribution of random variables \citep{chen2022robust}, which, once correctly identified, can accurately reflect the uncertain characteristics of random variables. However, obtaining the precise probability distribution information is closely related to the number of data samples as well as the accuracy of the  prediction methods. This increases the complexity of calculations and reduces confidence \citep{firouzmakan2019comprehensive}. Contrary to fuzzy and stochastic programming, RO defines the uncertain parameters in a problem as an uncertain set with bounds or scenarios, and makes decisions based on the average, or the worst case. Accordingly, after realizing any uncertain parameters, the optimal solution is achieved with good generalization, and does not have the mentioned complexity \citep{lu2020robust}.

RO considers a number of scenarios as possible outcomes, and assigns a probability of occurrence to each outcome. The objective of RO is to obtain a solution which is nearly optimal and feasible under each scenario \citep{4}. In this approach, one must still use historical data to find the probability of occurrence of each scenario, but this problem is solved by considering the possibility of solution infeasibility.

Stochastic programming optimizes the expected value of objective function over the set of all possibilities. On the contrary, robust optimization generates a solution which is feasible for all scenarios, and considers the infeasibility as a cost in the model \citep{shabani2014value, bertsimas2003robust}. In other words, although finding the correct probability of occurrence of the scenarios is still a challenge for robust optimization, it is dealt with by assuming the cost of infeasibility in the model. In stochastic models, failing to predict the right value of parameters leads to solution infeasibility, elevated costs, and futility of the decision support system. RO, however, intrinsically considers infeasibility with an assigned cost. In other words, RO facilitates the probability prediction challenge by delivering a single optimal solution that is feasible under all scenarios \citep{ahmadvand2022robust}

Although there exist numerous studies in the area of supplier selection, only a few papers have considered both green and non-green criteria at the same time. Also, to the best of our knowledge, this paper is the first study that employs robust optimization to handle uncertainties of the green supplier selection (GSS) problem and cap-and-trade mechanism for carbon emission. Based on the literature, RO is one of the best options to solve various problems in supply chain. In other words, as confirmed by the literature, RO is an efficient approach to deal with uncertainty. However, this approach has not been applied to the green supplier selection problem.

In other words, in this study, we present a GSS model which is embedded in a closed-loop supply chain framework. The model regards both green and non-green measurements in the presence of distinct quantitative and qualitative factors. Moreover, it is assumed that the market operates under the cap-and-trade regulations imposed by the government. To deal with the described problem, an RO approach is employed. It is shown that the proposed RO is able to generate solutions that are close to real cases.

This study aims to address the following research questions:
\begin{itemize}
    \item the impact of cap-and-trade, as an environmental mandate, on the supply chain;
    \item finding an effective approach for controlling the intrinsic uncertainty of the model's parameters;
    \item selecting the best supplier among all candidates in a supply chain considering economic and environmental criteria; and,
    \item determining clear cause-and-effect relationships between various factors to help government entities with their decision making efforts.
\end{itemize}

The structure of the remainder of the paper is as follows. Section \ref{Literature review} summarizes the relevant research. Section \ref{Mathematical model formulation} outlines the RO model of the considered GSS problem. Section \ref{Computational results} explains the computational experiments of the proposed model and discusses the results. Conclusions and future study directions are presented in section \ref{Conclusion}.

\section{Literature review}\label{Literature review}

The goal of green supply chain management is reduction of the harmful effects of the supply chain's activities on the environment. In this regard, firms have to identify the most effective measures for evaluating the environmental performance of their suppliers. Research in green supplier selection, in comparison with traditional supplier selection (which generally deals with the firm's profit and quality of the products) is limited. Hence, the common environmental criteria used in the literature of green supplier selection can be identified. Some of the most used environmental criteria in the literature include:

\begin{itemize}
    \item Environmental management system: this criteria defines the suppliers' policies for making the production process environmentally friendly; one example for such policies is the ISO 14001 certificate. \citep{8handfield2002,8hsu2009,8lee2009,8awasthi2010,8bai2010,8kuo2010,8mafakheri2011,8yeh2011,8amin2012, 8govindan2013,8tseng2013,10hu2015,10hashemi2015,8,10qin2017,10kumar2017,10gupta2017};
    
    \item Pollution production: this criteria measures the amount of pollution that a manufacturer creates \citep{8amin2012,8govindan2013,10hashemi2015,10hu2015,10kannan2015,10huang2016,8,10qin2017,10luthra2017,10kumar2017};
    \item Recyclability: refers to the capability of suppliers in using recycled material in their manufacturing process \citep{8yeh2011,8amin2012,10hu2015,10kannan2015,10govindan2016};
    \item Green product: this measure evaluates the ability of suppliers in using green technology as well as environmentally friendly material \citep{8handfield2002,8lee2009,8amin2012,8tseng2013};
    \item Product toxicity: assesses the suppliers regarding the level of toxic substance used in their products \citep{10hu2015,10kannan2015,8,10gupta2017}.
\end{itemize}

Table \ref{tab:literature review on GSS} summarizes the findings of the mentioned papers. One of the challenges of the supplier selection problem is that firms tend to maximize their profit while trying to have a great performance from the environmental point of view. Therefore, companies are required to identify the non-green criteria affecting the firm's profit as well as green criteria. Different studies have been completed to identify the most important criteria in traditional supplier selection problem. \cite{8dickson1966}, \cite{35lehmann1974}, \cite{35weber1991}, and \cite{8cheraghi2004} conducted research to identify non-green supplier evaluation measures. Those studies showed that cost, delivery performance, and quality are the three most used and important criteria in this particular problem area.

The next step of supplier selection and order allocation process is assessing the candidates regarding the criteria mentioned above. Researchers have applied different techniques to solve this problem. \cite{5chai2013} classified these techniques in three main groups: multi-criteria decision-making (MCDM), artificial intelligence (AI), and mathematical programming (MP). A list of the papers that apply these techniques to solve the supplier selection problem is presented in table \ref{tab:literature review on GSS}. MCDM is a framework that helps the decision makers find the best alternative between multiple options based on various criteria. To select the best option, MCDM sorts the alternatives based on their scores. AI refers to the science of making the computers able to work, learn, and think intelligently, like humans. MP is a useful approach to address the supply chain management problems in a clear way; constraints or equations can implement assumptions that make the problem realistic. In this paper, an MP is presented to model the problem.

Applied research invariably involves the development of methods to more closely align decisions to real-life situations. One such approach to augment applicability involves considering the supply chain as a closed-loop system. Firms in the supply chain can collect used or rejected products and deploy them again in the production system, thus contributing to reduced, environmentally harmful waste \citep{cao2020production}. Moreover, pollution control systems can also provide enhanced applicability for solutions. Cap-and-trade is an interesting mechanism that is becoming more popular in today’s consumer market. This mechanism is proved to be one of most effective approaches of GHG emission control \citep{xu2021channel}. Furthermore, based on \cite{yu2021reselling}, cap-and-trade is more effective than the carbon tax system on reducing GHG emissions because it produces a higher social welfare and generates more profit for the manufacturers. \cite{zhang2013multi} presented a multi-product production planning model under the C\&T system and proposed a profit-maximization model to achieve the firm's optimal policy. In their analysis, it was indicated that C\&T curbs the emission level better than carbon tax system. Another research supporting the effectiveness of C\&T is conducted by \cite{chen2020clean}, which further compared emission reduction effects of carbon tax and C\&T schemes. They showed that both mechanisms stimulate clean innovation, but C\&T is more efficient. They stated that government can control air pollution by assigning the proper carbon cap for manufacturers in a trade-off between environmental and economical objectives.

The emission of greenhouse gas was restricted by the Kyoto Protocol \citep{oberthur1999} for the first time. The Kyoto Protocol is an international treaty developed within the United Nations Framework Convention on Climate Change (UNFCCC). It was adopted in Kyoto, Japan in 1997 and came into effect in 2005. The Protocol proposes cap-and-trade, which is a flexible framework for reducing air pollution. 

The existing literature regards cap-and-trade as one of the most effective ways to control the carbon emission \citep{golpira2022robust}. \cite{zhang2013multi} considered cap-and-trade mechanism in their production planning and presented a profit maximization model to find the firm's best policy. \cite{gong2013optimal} designed a dynamic production model to scrutinize the impacts of creating a carbon trading market on production policy and found an emission trading policy as well as the optimal production plan under the cap-and-trade regulation. \cite{shen2014california} studied California's cap-and-trade scheme with the goal of implementing it as China's carbon reduction program. \cite{li2018game} studied the impact of cap-and-trade system on manufacturers' optimal operational decision and showed that customers' green preferences act as an incentive for greening the production technology. \cite{zhang2019evolutionary} considered two scenarios for carbon allowance prices: dynamic and static. They investigated the effects of cap-and-trade market on manufacturers' decisions under both scenarios and showed that upgrading the production technology is positively correlated with penalties imposed for extra emissions.

In the context of supply chain management, cap-and-trade is a relatively new topic. Thus, the need for more research in this area is deeply felt. To the best of the authors' knowledge cap-and-trade regulations have not been addressed in the research body related to green supplier selection.

\subsection{Research contributions}
Although there are numerous studies on the supplier selection problem, we know of only a small group of papers addressing both green and non-green factors in the process of supplier evaluation while considering a closed-loop structure with uncertain parameters. Furthermore, as the literature shows, there is no paper that considers RO for green supplier evaluation as well as the cap-and-trade mechanism. This paper is an attempt to fill these gaps by fostering the contributions presented below.
\begin{itemize}
    \item A multi-objective mathematical model for a closed-loop supplier selection and order allocation evaluating the candidates in terms of both environmental and economical criteria is presented. The developed model helps firms achieve a more environmentally friendly manufacturing system. Model realism is enhanced by developing an approach that considers two groups of conflicting criteria (green and non-green).
    \item The cap-and-trade mechanism, as a method to manage air pollution, is employed in the model. Analyses on the cap and market prices of carbon are performed to help firms and governments determine the values of the parameters to achieve better results. These analyses can subsequently translate to lower cost and carbon emission, as well as more environmentally friendly products. Moreover, based on the analysis conducted on the cap-and-trade approach, this mechanism is demonstrated as a proper approach for carbon emission reduction.
    \item The generalized model is solved by the RO approach to handle the uncertainty embedded in the problem. Sensitivity analysis has been conducted on two parameters to illustrate the trade-off between model robustness and solution robustness, and solution deviation. This can inform decision-makers of the best parameter values.
\end{itemize}

\begin{table}
\centering
\caption{Literature Review Summary}
\begin{adjustbox}{max width=\textwidth}
\renewcommand{\arraystretch}{1}
\newcommand*{\TitleParbox}[1]{\parbox[c]{10.75cm}{\raggedright #1}}%
\begin{tabular}{|l|ll|lll|lll|l|l|}
\multicolumn{1}{l}{}        &       & \multicolumn{1}{l}{}  &      &    & \multicolumn{1}{l}{}                                                &       &    & \multicolumn{1}{l}{}                                                   & \multicolumn{1}{l}{}         & \multicolumn{1}{l}{}                                                                                                                                                                                                                                                \\ 
\hline
{References} & \multicolumn{2}{l|}{Criteria} & \multicolumn{3}{l|}{\begin{tabular}[c]{@{}l@{}}Solving \\Approach\end{tabular}} & \multicolumn{3}{l|}{\begin{tabular}[c]{@{}l@{}}Uncertainty \\Approach\end{tabular}} & {Closed-loop} & {Findings}\\ 
\cline{2-9}
                            & Green & Non-green             & MCDM & MP & AI                                                                  & Fuzzy & SO & RO                                                                     &                              &                                                                                                                                                                                                                                                                     \\ 
\hline
\cite{8handfield2002}                           & \checkmark     & -                     & \checkmark    &    &                                                                     &       &    &                                                                        &                              & AHP as a decision support model is used to  find out the trade-offs between environmental dimensions.                                                                                                                                                               \\ 
\hline
\cite{8hsu2009}                           & \checkmark     & -                     & \checkmark    &    &                                                                     &       &    &                                                                        &                              & \begin{tabular}[c]{@{}l@{}}Hazardous substance management (HSM) is incorporated into supplier selection using ANP approach.\end{tabular}                                                                                                                          \\ 
\hline
\cite{8lee2009}                           & \checkmark     & \checkmark                     & \checkmark    &    &                                                                     & \checkmark     &    &                                                                        &                              & The Delfi method and fuzzy AHP are used to support  supplier selection decision.                                                                                                                                                                                    \\ 
\hline
\cite{8awasthi2010}                           & \checkmark     & -                     & \checkmark    &    &                                                                     & \checkmark     &    &                                                                        &                              & \begin{tabular}[c]{@{}l@{}}A fuzzy multi-criteria approach is presented that consists three steps: finding the best criteria, scoring\\ the suppliers using fuzzy TOPSIS, and sensitivity analysis of criteria weights on supliers evaluation.\end{tabular}       \\ 
\hline
\cite{8bai2010}                           & \checkmark     & \checkmark                     &      &    & \checkmark                                                                   &       &    &                                                                        &                              & \begin{tabular}[c]{@{}l@{}}A new supplier selection technique using rough set and grey system theory is presented.\end{tabular}                                                                                                                                   \\ 
\hline
\cite{8kuo2010}                           & \checkmark     & \checkmark                     &      &    & \checkmark                                                                   &       &    &                                                                        &                              & \begin{tabular}[c]{@{}l@{}}A green supplier selection model is generalized using a hybrid method called ANN-MADA.\end{tabular}                                                                                                                                    \\ 
\hline
\cite{8mafakheri2011}                           & \checkmark     & \checkmark                     & \checkmark    & \checkmark  &                                                                     &       &    &                                                                        &                              & \begin{tabular}[c]{@{}l@{}}A two-stage multi-criteria dynamic programming approach for suppier \\selection and order allocation is proposed.\end{tabular}                                                                                                           \\ 
\hline
\cite{8yeh2011}                           & \checkmark     & \checkmark                     &      &    & \checkmark                                                                   &       &    &                                                                        &                              & \begin{tabular}[c]{@{}l@{}}A green partner selection model is solved by genetic algorithm to find the \\set of pareto optimal solutions.\end{tabular}                                                                                                               \\ 
\hline
\cite{8amin2012}                           & \checkmark     & \checkmark                     &      & \checkmark  &                                                                     & \checkmark     &    &                                                                        & \checkmark                            & A framework for supplier evaluation in a closed-loop supply chain is generalized.                                                                                                                                                                                   \\ 
\hline
\cite{8govindan2013}                          & \checkmark     & \checkmark                     & \checkmark    &    &                                                                     & \checkmark     &    &                                                                        &                              & \begin{tabular}[c]{@{}l@{}}Triple Bottom Line approach for supplier selection using a fuzzy multi-criteria \\model is developed.\end{tabular}                                                                                                                       \\ 
\hline
\cite{8tseng2013}                          & \checkmark     & \checkmark                     &      &    & \checkmark                                                                   &       &    &                                                                        &                              & \begin{tabular}[c]{@{}l@{}}Environmental and non-environmental~criteria for selecting the best partners \\are identified by evauating the weight of criteria and using grey relational analysis.\end{tabular}                                                       \\ 
\hline
\cite{10hu2015}                          & \checkmark     & \checkmark                     & \checkmark    &    &                                                                     &       &    &                                                                        &                              & \begin{tabular}[c]{@{}l@{}}A novel evaluation system of green supplier selection under the mode of low \\carbon economy is investigated.\end{tabular}                                                                                                               \\ 
\hline
\cite{10hashemi2015}                          & \checkmark     & \checkmark                     & \checkmark    &    &                                                                     &       &    &                                                                        &                              & \begin{tabular}[c]{@{}l@{}}A comprehensive green supplier selection model using ANP and Grey relational \\analysis is proposed.\end{tabular}                                                                                                                        \\ 
\hline
\cite{10kannan2015}                          & \checkmark     & \checkmark                     & \checkmark    &    &                                                                     & \checkmark     &    &                                                                        &                              & \begin{tabular}[c]{@{}l@{}}Fuzzy Axiomatic Design is proposed to to select the best green supplier for a plastic \\manufacturing company.\end{tabular}                                                                                                              \\ 
\hline
\cite{10huang2016}                          & \checkmark     & \checkmark                     &      & \checkmark  & \checkmark                                                                   &       &    &                                                                        &                              & \begin{tabular}[c]{@{}l@{}}A game-theoretic model is presented in order to investigate the impacts of supplier\\ selection, transportation mode selection, the product line design, and pricing strategies\\ on profits and greenhouse gases emissions.\end{tabular}  \\ 
\hline
\cite{8}                          & \checkmark     & \checkmark                     & \checkmark    &    &                                                                     &       &    &                                                                        &                              & Best worst method is used to find the best suppliers among the qualified suppliers.                                                                                                                                                                                 \\ 
\hline
\cite{10luthra2017}                          & \checkmark     & \checkmark                     & \checkmark    &    &                                                                     &       &    &                                                                        &                              & \begin{tabular}[c]{@{}l@{}}An integrated approach of AHP, VIKOR, and multi-criteria optimization is developed \\to solve evaluate the sustainable supplier selection.\end{tabular}                                                                                  \\ 
\hline
\cite{10qin2017}                          & \checkmark     & -                     & \checkmark    &    &                                                                     & \checkmark     &    &                                                                        &                              & \begin{tabular}[c]{@{}l@{}}TOMID approach to solve a green supplier selection problem considering interval \\type-2 fuzzy sets is extended.\end{tabular}                                                                                                            \\ 
\hline
\cite{10kumar2017}                          & \checkmark     & \checkmark                     & \checkmark    &    &                                                                     & \checkmark     &    &                                                                        &                              & \begin{tabular}[c]{@{}l@{}}Suppliers’ performance is evaluated based on Green Practices using the \\fuzzy-extended Elimination and Choice Expressing Reality approach.\end{tabular}                                                                                 \\ 
\hline
\cite{10gupta2017}                          & \checkmark     & -                     & \checkmark    &    &                                                                     & \checkmark     &    &                                                                        &                              & \begin{tabular}[c]{@{}l@{}}Supplier evaluation is done usnig a three phase methodology including three phases: \\identifying green criteria, ranking determined criteria using a novel best worst \\method, and ranking suppliers using fuzzy TOPSIS.\end{tabular}  \\ 
\hline
\cite{36}                          & \checkmark     & \checkmark                     &      & \checkmark  &                                                                     &       & \checkmark  &                                                                        &                              & \begin{tabular}[c]{@{}l@{}}A sustainable and efficient supply chain is designed using a multi-objective \\MINLP model for supplier selection and order allocation with stochastic demand.\end{tabular}                                                              \\ 
\hline
This paper                          & \checkmark     & \checkmark                     &      & \checkmark  &                                                                     &       &    & \checkmark                                                                      & \checkmark                            &                                                                                                                                                                                                                                                                     \\
\hline
\end{tabular}
\end{adjustbox}
\end{table}

\section{Mathematical model formulation}\label{Mathematical model formulation}
\subsection{Robust optimization}\label{Robust optimization}

An interesting contribution that makes the results of the supplier selection problem valuable is considering the intrinsic uncertainty of the parameters; wrong  estimation of parameter values may lead to colossal losses in an uncertain environment \citep{wang2021incentive}. Obviously, this requires an effective approach for handling such uncertainty. To the best of our knowledge, RO has not been applied to solve the problems in the literature of green supplier evaluation. In this paper, RO is employed because it handles infeasibility, and is confirmed to be an efficient approach to deal with data uncertainty by treating uncertainty as a deterministic trait, while not limiting uncertain parameter values to point estimates \citep{jeyakumar2014support}. Accordingly, to decrease its possibility, RO assigns a cost to infeasibility. Other popular approaches are capable of controlling uncertainty, but in case of solution infeasibility, which may tremendously increase the costs, other approaches are required. In other words, the expected cost of infeasibility can be significant and must be taken into consideration.

Unlike deterministic methods that feature a specific value for the parameters, RO considers various parameter values and obtains solutions under all possible scenarios. \cite{22mulvey1995} introduced RO as a new approach for uncertain problems. RO performs a trade-off between two kinds of robustness; namely, model robustness (i.e., the obtained solution yields the least infeasibility in all scenarios), and solution robustness(i.e., the solution is as close to the optimal value as possible in all scenarios). Simply put, the goal of RO is to achieve a robust solution which ensures all the real scenarios are nearly optimal and feasible. 

\cite{22mulvey1995} assumed $x$ and $y$ as design variables (in which their optimal value is not dependent upon the real value of the uncertain parameters), and control variables (in which their optimal value depends on the value of the uncertain parameters). Correspondingly, there are two groups of constraints: control constraints, and structural constraints. The former are those that include noisy parameters, while the latter feature the constraints with no uncertain parameters. The basic structure of the RO model is presented below.
\begin{align}
& Min \; \;c^{T}x+d^{T}y \label{eq:1 obj of robust examplanation} \\
& Ax=b \label{eq:2 first constraint of robust examplanation} \\
& Bx+Cy=e \label{eq:3 second constraint of robust examplanation} \\
& x,y\geq0 \label{eq:4 positive variables in robust examplanation}
\end{align}
where $A$ and $b$ are deterministic parameters, and $B$ and $C$ are uncertain parameters. In other word, equations (\ref{eq:2 first constraint of robust examplanation}) and (\ref{eq:3 second constraint of robust examplanation}) convey the structural and control constraints, respectively. There are a finite set of scenarios $s$ defined for the RO problem with the probability $Pr_{s}$ for each scenario. Then, the model formulation transforms to the following: 
\begin{align}
& Min \; \;\sigma(x,y_{1},...,y_{s})+\omega\rho(\delta_{1},...,\delta_{s})\label{eq:5 objective of second robust examplanation}\\
& Ax=b \label{eq:6 deterministic constraint in second robust examplanation}\\
& B_{s}x+C_{s}y+\delta_{s}=e_{s}\label{eq:7 uncertain constraint in second robust examplanation}\\
& x,y_{s}\geq0\label{eq:8 positive uncertain variables in robust explanation}
\end{align}
in which $\delta_{s}$ indicates the permitted infeasibility in constraint (\ref{eq:7 uncertain constraint in second robust examplanation}) under scenario $s$. The objective function includes two parts. The first part represents the solution robustness, and the second part is a measure of model robustness. As the value of the penalty of infeasibility, $\omega$ helps to make a trade-off between model robustness and solution robustness. The first part of the presented objective function can be reformulated as follows:
\begin{equation}\label{eq:9 quadratic solution robustness}
\sigma(x,y_{1},...,y_{s})=\sum_{s}Pr_{s}\xi_{s}+\lambda \sum_{s}Pr_{s}\left(\xi_{s}-\sum_{s'}Pr_{s'}\xi_{s'}\right)^{2}
\end{equation}
where $\lambda$ is the penalty for solution variance and takes positive values. $\xi_{s}$ is the value of the objective function in each scenario. In this equation, the first part shows the average value of the objective in different scenarios, and the second part indicates the objective variance. Furthermore, equation (\ref{eq:9 quadratic solution robustness}) can be written in the following form: 
\begin{equation}\label{eq:10 modulus solution robustness}
\sigma(x,y_{1},...,y_{s})=\sum_{s}Pr_{s}\xi_{s}+\lambda \sum_{s}Pr_{s} \mid \xi_{s}-\sum_{s'}Pr_{s'}\xi_{s'}\mid
\end{equation}

\cite{4yu2000} state that, due to non-linearity, this equation is very complicated and time-consuming to solve. Instead, they propose the following formulation to linearize the equation.
\begin{align}
&\sigma(x,y_{1},...,y_{s})=\sum_{s}Pr_{s}\xi_{s}+\lambda \sum_{s}Pr_{s}\left[\left(\xi_{s}-\sum_{s'}Pr_{s'}\xi_{s'}\right)+2\theta_{s}\right]  \label{eq:11 linearized solution robustness}\\
&\left(\xi_{s}-\sum_{s'}Pr_{s'}\xi_{s'}\right)+\theta_{s}\geq0 \qquad\qquad \forall s\label{eq:12 auxilary constraint for solution robustness}\\
&\theta_{s}\geq0 \qquad\qquad \forall s\label{eq:13 theta positive}
\end{align}

The variable $\theta_{s}$ is used to prevent $\xi_{s}-\sum_{s'}P_{s'}\xi_{s'}$ from being negative. According to the above new formulation, the objective function can be written as follows.
\begin{equation}\label{eq:14 general obj in robust models}
Min\;\;\sum_{s}Pr_{s}\xi_{s}+\lambda \sum_{s}Pr_{s}\left[\left(\xi_{s}-\sum_{s'}Pr_{s'}\xi_{s'}\right)+2\theta_{s}\right]+\omega\sum_{s}Pr_{s}\delta_{s}
\end{equation}
where $\omega$ is defined as the weight of trade-off between model robustness and solution robustness. $\omega$ is a number without unit, which indicates the penalty assigned to infeasibility; if $\omega$ increases, the probability of obtaining an infeasible solution decreases.

\subsection{Model description}\label{Model description}
In this section, a multi-objective, multi-product, multi-supplier, and multi-period robust model for green supplier selection and order allocation in a closed-loop supply chain is presented. Two groups of criteria are presented for supplier evaluation: green and non-green. The non-green criteria are cost, quality, and delivery time, while green criteria are divided into two categories: quantitative (carbon emission amount), and qualitative (environmental management system level, recyclability level, green product level, product toxicity level). 

In this closed-loop supply chain, the manufacturer sells its products to the customer, but some products may be rejected due to their low quality or other problems. Moreover, the manufacturer tends to use second-hand products in the production process to benefit from their low costs. Second-hand products must be collected from the customers. A percentage of the collected and rejected products are usable. The reminder are disposed. Therefore, the input material and component parts include product purchased from suppliers, reusable parts or material collected from the final consumer market, and usable parts from rejected products (see figure \ref{Fig: Closed-loop supply chain structure}). Other assumptions of the model are presented below.

\textbf{Assumptions}
    \begin{itemize}
    \item Shortage is allowed and will be back-ordered. A penalty will be charged for each unit of shortage.
    \item The supply chain operates under cap-and-trade regulations. According to this mechanism, each manufacturer is allowed to emit a predetermined amount of carbon. Any extra quantity of carbon emission must be purchased from other manufacturers in a trade market. Conversely, manufacturers can sell their unused carbon emission allowance to generate profit. Note that a carbon tax system allocates a price to carbon emission and puts it on the market to find the emission reduction mechanism. On the other hand, cap-and-trade assigns an emission level to the manufacturers and allows the market to find the allowance price. So, in cap-and-trade mechanism. In other words, buyers and sellers offer their desired prices, which can be different based on their profit level, and the manufacturer considers the maximum offer for selling its allowance and the minimum offer for buying extra allowance.
    \item Uncertain parameters include demand, cost of products purchased from suppliers,  percentage of the returned products, percentage of the collected second-hand products, percentage of reusable collected products, percentage of usable rejected products, and delivery delay.
    \item Carbon emission is a function of $CO_{2}$ released in the production process as well as $CO_{2}$ released during material and product transportation. Purchased products from a supplier can be transported by either the manufacturer's trucks (buyer) or the supplier's trucks. The manufacturer is responsible for the carbon emission during the transportation process only if the manufacturer's trucks are employed. Using the manufacturer's truck for transportation is cheaper but increases the amount of carbon released by the manufacturer.
    \item Different kinds of trucks exist for transporting the procured materials and products according to the weight or volume of the loads. 
     \item A constant interest rate is applied to all of the costs through time.
    \end{itemize}
    
\begin{figure}[h]
\includegraphics[width=\textwidth]{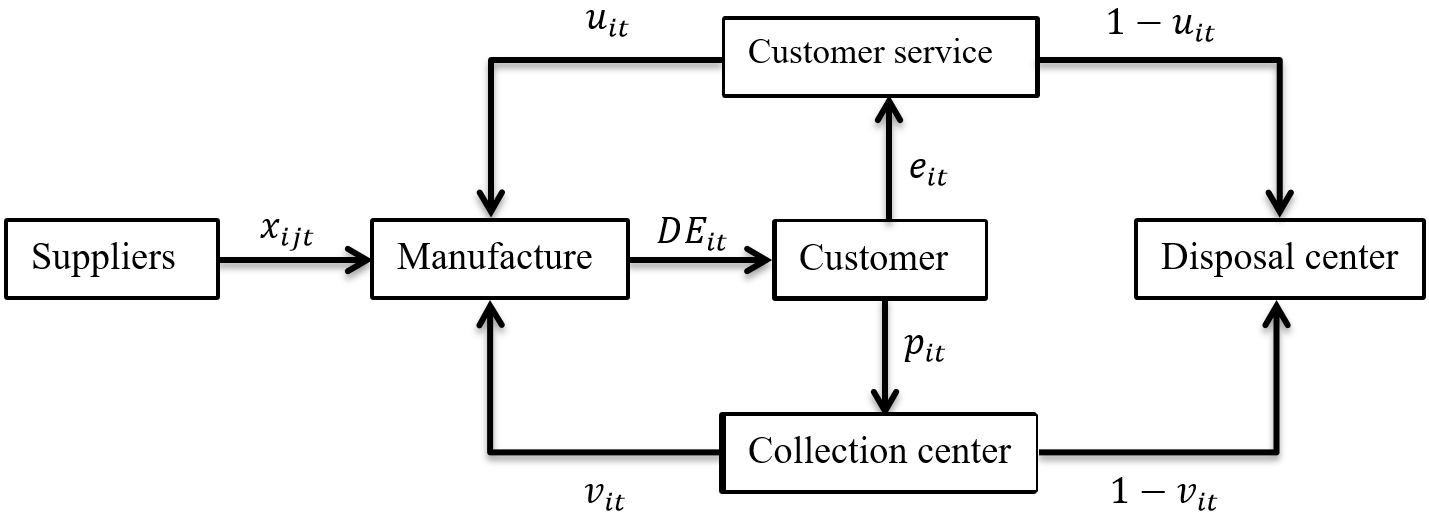}
\caption{Closed-loop supply chain structure}
\label{Fig: Closed-loop supply chain structure}
\end{figure}

\textbf{Indices}

\noindent $i$: products $(i=1,\cdots,I)$

\noindent $j$: suppliers $(j=1,\cdots,J)$

\noindent $t$: periods $(t=1,\cdots,T)$

\noindent $k$: vehicles $(k=1,\cdots,K)$

\noindent $n$: echelons in supply chain $(n=1,2)$ ($n=1$ and $n=2$ show the supplier and buyer, respectively).

\noindent $m$: set of prices offered by buyers or sellers at the carbon emission trading market

\noindent $s$: scenarios $(s=1,\cdots,S)$

\textbf{Parameters}

\noindent $Pr_{s}$: probability of scenario $s$

\noindent $C_{ijt}^{s}$: cost of product $i$ purchased from supplier $j$ in period $t$ under scenario $s$

\noindent $ir$: interest rate

\noindent $h_{it}$: holding cost of product $i$ in period $t$

\noindent $f_{it}$: backorder cost of product $i$ in period $t$

\noindent $L_{jt}^{s}$: number of days in which products purchased from supplier $j$ in period $t$ under scenario $s$ received after the specified time

\noindent $G_{ijt}$: penalty for each day of late delivery of product $i$ purchased from supplier $j$ in period $t$ 

\noindent $e_{it}^{s}$: percentage of product $i$ rejected in period $t$ under scenario $s$ from the customer due to low quality

\noindent $u_{it}^{s}$: percentage of usable parts after disassembly of product $i$ rejected by customer in period $t$ under scenario $s$

\noindent $p_{it}^{s}$: percentage of used product $i$ collected in period $t$ under scenario $s$ from the customer 

\noindent $v_{it}^{s}$: percentage of reusable parts after disassembly of product $i$ collected as used product in period $t$ under scenario $s$

\noindent $O_{ijt}$: manufacturer's loss caused by rejected product $i$ purchased from supplier $j$ in period $t$

\noindent $SP_{t}^{m}$: price offered by seller $m$ in the carbon emission allowance market in period $t$

\noindent $BP_{t}^{m}$: price offered by buyer $m$ in the carbon emission allowance market in period $t$

\noindent $DC_{it}$: disassembly cost for product $i$ in period $t$

\noindent $RC_{it}$: remanufacturing cost of product $i$ in period $t$

\noindent $DP_{it}$: disposal cost of product $i$ in period $t$

\noindent $TC_{jtkn}$: transportation cost of products purchased from supplier $j$ in period $t$ by truck type $k$ belonging to the $n$th echelon of supply chain

\noindent $DE_{it}^{s}$: customers' demand for product $i$ in period $t$ under scenario $s$

\noindent $d_{j}$: distance of supplier $j$ from the producer

\noindent $CET_{jtkn}$: amount of carbon emitted per kilometer of transportation from supplier $j$ in period $t$ by truck type $k$ belonging to $n$th echelon of supply chain

\noindent $CEP_{ijt}$: amount of carbon emitted during production using product $i$ purchased from supplier $j$ in period $t$

\noindent $CER_{it}$: amount of carbon emitted during remanufacturing or recycling product $i$ in period $t$

\noindent $EM_{ijt}$: environmental management system score assigned to supplier $j$ from which product $i$ in period $t$ is purchased. This score is determined by experts at the manufacturer's firm evaluating the performance of the suppliers regarding environmental policies(e.g. possession of an ISO 14001 certificate). 

\noindent $GP_{ijt}$: green product score of supplier $j$ from which product $i$ in period $t$ is purchased. This is the manufacturer's strategy to purchase products with minimum environmental effect during its life cycle.

\noindent $RE_{ijt}$: recyclability score of product $i$ purchased from supplier $j$ in period $t$

\noindent $PT_{ijt}$: toxicity score of product $i$ purchased from supplier $j$ in period $t$

\noindent $CAP_{t}$: amount of carbon emission determined by the government as an upper bound for the manufacturer in period $t$

\noindent $M_{k}$: breaking point of loads for each kind of truck used for transportation

\textbf{Decision variables}

\noindent $x_{ijt}$: number of product $i$ purchased from supplier $j$ in period $t$

\noindent $r_{it}$: number of remaining product $i$ at the end of period $t$

\noindent $b_{it}$: number of back-ordered product $i$ at the end of period $t$

\noindent $\alpha_{t}$: allowance of carbon emission the manufacturer needs to buy from the market in period $t$

\noindent $\beta_{t}$: allowance of carbon emission the manufacturer wants to sell to the other manufacturers in period $t$

\noindent $\tau_{t}$: the best price for the manufacturer to buy carbon allowance in period $t$

\noindent $\phi_{t}$: the best price for the manufacturer to sell carbon allowance in period $t$

\noindent $q_{jtkn}$: binary variable, equals 1 if the buyer places an order to supplier $j$ in period $t$ transported by truck $k$ belonging to $n$th echelon; 0 otherwise

\noindent $W_{jtkn}$: auxiliary continuous variable to determine the amount of load from supplier $j$ in period $t$ transported by truck type $k$ belonging to $n$th echelon

\noindent $\delta_{ijt}^{s+}$: under-fulfillment (shortage) of product $i$ purchased from supplier $j$ in period $t$ 

\noindent $\delta_{ijt}^{s-}$: over-fulfillment (storage) of product $i$ purchased from supplier $j$ in period $t$ 

According to the criteria introduced, three objectives are defined. The first objective, $\xi_{1}^{s}$, is cost, which is a measure of non-green criteria. this cost includes uncertain parameters and is formulated as:
\begin{equation}\label{eq:15 first obj}
\begin{aligned}
\underset{s \in \{1,2,\ldots,S\}}{\xi_{1}^{s}} &=\sum_{j}\sum_{t}L_{jt}^{s}\sum_{i}G_{ijt}x_{ijt}+\sum_{i}\sum_{t}e_{it}^{s}\sum_{j}O_{ijt}x_{ijt}+\sum_{i}\sum_{j}\sum_{t}x_{ijt}C_{ijt}^{s}\\ &+\sum_{i}\sum_{t}r_{it}h_{it}+\sum_{i}\sum_{t}b_{it}f_{it}+\sum_{i}\sum_{t}e_{it}^{s}\sum_{j}x_{ijt}DC_{it}\\
&+\sum_{i}\sum_{t}p_{it}^{s}\sum_{j}x_{ijt}DC_{it}+\sum_{i}\sum_{t}e_{it}^{s}u_{it}^{s}RC_{it}\sum_{j}x_{ijt}\\
&+\sum_{i}\sum_{t}p_{it}^{s}v_{it}^{s}RC_{it}\sum_{j}x_{ijt}+\sum_{i}\sum_{t}e_{it}^{s}(1-u_{it}^{s})DP_{it}\sum_{j}x_{ijt}\\
&+\sum_{i}\sum_{t}p_{it}^{s}(1-v_{it}^{s})DP_{it}\sum_{j}x_{ijt}+\sum_{j}\sum_{t}\sum_{k}\sum_{n}TC_{jtkn}q_{jtkn}\\
&-\sum_{t}\alpha_{t}\tau_{t}+\sum_{t}\beta_{t}\phi_{t}
\end{aligned}
\end{equation}

In this equation, the first term calculates the penalty for rejected items. It is an indication of the quality of the products. The second term is the cost of delay. The remaining terms calculate other features including purchasing cost, holding cost, back-order cost, disassembly cost of rejected and collected products, remanufacturing cost for the usable parts obtained from collected and rejected products, disposal cost of the unusable parts, and profit/cost of selling/buying carbon emission allowance. Of note, three different kinds of criteria (cost, quality, and delivery performance) in this objective, all of which are combined in the cost function $\xi_{1}^{s}$.

The second objective is the amount of carbon emission, which is a measure of quantitative green criteria. Due to parameter uncertainty $\xi_{2}^{s}$ is written as:

\begin{equation}\label{eq:16 second objective}
\begin{aligned}
\underset{s \in \{1,2,\ldots,S\}}{\xi_{2}^{s}}
& =\sum_{j} d_{j}\sum_{k}\sum_{t}CET_{jtk2}q_{jtk2}+\sum_{i}\sum_{j}\sum_{t}x_{ijt}CEP_{ijt}\\
&+\sum_{i}\sum_{t}(e_{it}^{s}u_{it}^{s}+p_{it}^{s}v_{it}^{s})CER_{it}\sum_{j}x_{ijt}
\end{aligned}
\end{equation}

The first term of equation (\ref{eq:16 second objective}) calculates the amount of carbon released during transportation of purchased products from supplier $j$ by trucks belonging to the manufacturer. The second term measures the carbon emission during the production phase, and the third term calculates the carbon released during the remanufacturing or recycling process fed by components and parts obtained from disassembling the collected and rejected products.

The third objective is a measure of qualitative green criteria. This objective function does not include uncertain parameters, and is defined as:

\begin{equation}\label{eq:17 third objective}
\begin{aligned}
z_{3} &=\sum_{i}\sum_{j}\sum_{t}LC_{ijt}x_{ijt}+\sum_{i}\sum_{j}\sum_{t}ER_{ijt}x_{ijt}+\sum_{i}\sum_{j}\sum_{t}RE_{ijt}x_{ijt}\\
&+\sum_{i}\sum_{j}\sum_{t}TS_{ijt}x_{ijt}
\end{aligned}
\end{equation}

The robust model of the described problem is therefore formulated as:
 
\begingroup
\allowdisplaybreaks
\begin{align}
\begin{split}
    & Min \; z_{1} =\sum_{s}Pr_{s}\xi_{1}^{s} + \lambda_{1} \sum_{s}Pr_{s} \left[\left(\xi_{1}^{s}-\sum_{s'}Pr_{s'}\xi_{1}^{s'}\right)+2\theta_{1}^{s}\right] \label{eq:18 first robust model obj}\\
    &  \qquad\qquad\qquad +\omega \sum_{s}Pr_{s} \left( \sum_{i}\sum_{j}\sum_{t} \left(\delta_{ijt}^{s^{-}}+\delta_{ijt}^{s^{+}}\right)\right)
\end{split}\\[2ex]
\begin{split}
    & Min \;z_{2} = \sum_{s}Pr_{s}\xi_{2}^{s}+\lambda_{2} \sum_{s}Pr_{s}\left[\left(\xi_{2}^{s}-\sum_{s}Pr_{s'}\xi_{2}^{s'}\right)+2\theta_{2}^{s}\right]\label{eq:19 second robust model obj}\\
    &  \qquad\qquad\qquad +\omega\sum_{s}Pr_{s}\left( \sum_{i}\sum_{j}\sum_{t}\left(\delta_{ijt}^{s^{-}}+\delta_{ijt}^{s^{+}}\right)\right)
\end{split}\\[2ex]
& Max\; z_{3}\label{eq:20 third robust model obj}\\[2ex]
& \left(\xi_{1}^{s}-\sum_{s'}Pr_{s'}\xi_{1}^{s'}\right)+\theta_{1}^{s}\geq 0 \qquad\qquad\qquad\qquad\qquad\qquad\qquad\qquad \forall s\label{eq:21 auxilary constraint for first obj}\\[2ex]
& \left(\xi_{2}^{s}-\sum_{s'}Pr_{s'}\xi_{2}^{s'}\right)+\theta_{2}^{s}\geq 0 \qquad\qquad\qquad\qquad\qquad\qquad\qquad\qquad \forall s \label{eq:22 auxilary constraint for second obj}\\[2ex]
\begin{split}
    & \sum_{j}d_{j}\sum_{k}CET_{jtk2}q_{jtk2}+\sum_{i}\sum_{j}x_{ijt}CEP_{ijt}\\ \label{eq:23 cap constraint}
     & \qquad\qquad
     +\sum_{i}(e_{it}^{s}u_{it}^{s}+p_{it}^{s}v_{it}^{s})CER_{it}\sum_{j}x_{ijt} \leq CAP_{t}-\alpha_{t}+\beta_{t} \qquad\qquad  \forall t,s
\end{split}\\[2ex]
&\tau_{t}= \max_{\forall m} BP_{t}^{m} \qquad\qquad\qquad\qquad\qquad\qquad\qquad\qquad\qquad\qquad\qquad  \forall t\label{eq:24 maximum price detection}\\[2ex]
& \phi_{t}= \min_{\forall m} SP_{t}^{m} \qquad\qquad\qquad\qquad\qquad\qquad\qquad\qquad\qquad\qquad\qquad  \forall t\label{eq:25 minimum price detection}\\[2ex]
\begin{split}
    &u_{it}^{s}e_{it}^{s}\sum_{j}x_{ijt}+v_{it}^{s}p_{it}^{s}\sum_{j}x_{ijt}+\sum_{j}x_{ijt}+b_{it}+r_{it-1}+\sum_{j}\delta_{ijt}^{s-}\label{eq:26 inventory balance}\\
    &  \qquad\qquad =DE_{it}^{s}+r_{it}+b_{it-1}+\sum_{j}\delta_{ijt}^{s+} \qquad\qquad\qquad\qquad\qquad\forall i,t,s
\end{split}\\[2ex]
& \sum_{i}x_{ijt}=\sum_{k}M_{k}\sum_{n}W_{jtkn} \qquad\qquad\qquad\qquad\qquad\qquad\qquad\qquad \forall j,t\label{eq:27 category calculation1}\\[2ex]
& W_{jt1n}\leq q_{jt1n} \qquad\qquad\qquad\qquad\qquad\qquad\qquad\qquad\qquad\qquad\qquad \forall j,t,n\label{eq:28 category calculation2}\\[2ex]
& W_{jt(k+1)n}\leq q_{jt(k+1)n}+q_{jtkn} \qquad\qquad\qquad\qquad\qquad \forall j,t,n;(k=1,\cdots,(K-1)) \label{eq:29 category calculation3}\\[2ex]
& W_{jtKn}\leq q_{jt(K-1)n} \qquad\qquad\qquad\qquad\qquad\qquad\qquad\qquad\qquad\qquad \forall j,t,n\label{eq:30 category calculation4}\\[2ex]
& \sum_{k}q_{jtkn}=1 \qquad\qquad\qquad\qquad\qquad\qquad\qquad\qquad\qquad\qquad\qquad \forall j,t,n\label{eq:31 category calculation5}\\[2ex]
& \sum_{k}W_{jtkn}=1 \qquad\qquad\qquad\qquad\qquad\qquad\qquad\qquad\qquad\qquad\qquad \forall j,t,n\label{eq:32 category calculation6}\\[2ex]
& 0\leq W_{jtkn}\leq 1 \qquad\qquad\qquad\qquad\qquad\qquad\qquad\qquad\qquad\qquad\qquad \forall j,k,t,n\label{eq:33 category calculation7}\\[2ex]
& C_{ij(t+1)}^{s}=C_{ijt}^{s}(1+ir) \qquad\qquad\qquad\qquad\qquad\qquad\qquad\qquad\qquad \forall i,j,t,s\label{eq:34 product cost inflaion}\\[2ex]
& h_{i(t+1)}=h_{it}(1+ir) \qquad\qquad\qquad\qquad\qquad\qquad\qquad\qquad\qquad\qquad \forall i,t\label{eq:35 holding cost inflaion}\\[2ex]
& f_{i(t+1)}=f_{it}(1+ir) \qquad\qquad\qquad\qquad\qquad\qquad\qquad\qquad\qquad\qquad \forall i,t\label{eq:36 backorder cost inflaion}\\[2ex]
& G_{ij(t+1)}=G_{ijt}(1+ir) \qquad\qquad\qquad\qquad\qquad\qquad\qquad\qquad\qquad \forall i,j,t\label{eq:37 late delivery penalty inflaion}\\[2ex]
& O_{ij(t+1)}=O_{ijt}(1+ir) \qquad\qquad\qquad\qquad\qquad\qquad\qquad\qquad\qquad \forall i,j,t\label{eq:38 rejection penalty inflaion}\\[2ex]
& DC_{i(t+1)}=DC_{it}(1+ir) \qquad\qquad\qquad\qquad\qquad\qquad\qquad\qquad\qquad \forall i,t\label{eq:39 disassembly cost inflaion}\\[2ex]
& RC_{i(t+1)}=RC_{it}(1+ir) \qquad\qquad\qquad\qquad\qquad\qquad\qquad\qquad\qquad \forall i,t\label{eq:40 remanufacturing cost inflaion}\\[2ex]
& DP_{i(t+1)}=DP_{it}(1+ir) \qquad\qquad\qquad\qquad\qquad\qquad\qquad\qquad\qquad \forall i,t\label{eq:41 disposal cost inflaion}\\[2ex]
& TC_{j(t+1)kn}=TC_{jtkn}(1+ir) \qquad\qquad\qquad\qquad\qquad\qquad\qquad\qquad \forall j,t,k\label{eq:42 transportation cost inflaion}
\end{align}
\endgroup

Equation (\ref{eq:18 first robust model obj}) aims to minimize the costs while trying to obtain the highest profit in the carbon emission trading market. Due to defining different scenarios, this objective is written as a robust objective to obtain the closest solution to all scenarios while keeping infeasibility at the lowest level. Equation (\ref{eq:19 second robust model obj}) is another robust objective aiming to minimize the carbon emission. Note that equations \eqref{eq:18 first robust model obj} and \eqref{eq:19 second robust model obj} go in opposite directions. In other words, the first and second objectives restrict each other. Therefore, the optimal solution is found through a trade-off between these two objectives. The third objective maximizes quality by choosing the suppliers with the best qualitative performance. Suppliers' qualitative scores are defined by experts on a scale of $0-10$. 

Equations (\ref{eq:21 auxilary constraint for first obj}) and (\ref{eq:22 auxilary constraint for second obj}) transform the first and second objectives to linear functions. These equations ensure that the deviation of each scenario from the average objective value is a positive amount. Equation (\ref{eq:23 cap constraint}) adheres to the manufacturer's carbon emission allowance, yet makes it possible to buy or sell the allowance. Buying and selling price is defined in equations (\ref{eq:24 maximum price detection}) and (\ref{eq:25 minimum price detection}), where the maximum and minimum offered price is selected to sell and buy the allowance, respectively. Equation (\ref{eq:26 inventory balance}) balances the inventories. 

Equations (\ref{eq:27 category calculation1}) to (\ref{eq:33 category calculation7}) confirm that the purchased products from each supplier are assigned to a specific truck category based on their size. For instance, if three different kinds of trucks are classified by their load size, and the order from a specific supplier fits on the first category, the variable $x_{ijt}$ will be a linear combination of first and second breaking points. Hence, $q_{jt1n}$ will be equal to 1 and variables $W_{jt1n}$ and $W_{jt2n}$ will be positive. Equations (\ref{eq:34 product cost inflaion})-(\ref{eq:42 transportation cost inflaion}) refer to the interest rate applied to the prices.

\subsection{Solution procedure}\label{Solution procedure}
To solve the proposed RO model a two-step procedure was followed. In the first step, each of the three objectives are independently solved , i.e., the model is solved three times; each time, one of the objectives is considered and the other two objectives are removed from the model. This way, the optimum value of each objective is found in the absence of the other objectives. Obviously, in multi-objective models the value of each objective (in the presence of the other objectives) is worse than or equal to their global optimal value. Thus, there is always a deviation between objective function values in multi-objective models and their optimal values; the goal is to minimize this deviation, which will be done in the second step. In the second step, the model is reformulated as a single objective mathematical model whose objective function is minimizing the normalized deviation of the objectives from their optimal values. Assume that the values of the first, second, and third objective are demonstrated as $z_{1}$, $z_{2}$, and $z_{3}$, respectively; and the optimal values are illustrated as $z_{1}^{*}$, $z_{2}^{*}$, and $z_{3}^{*}$. Then, the objective function is as follows.
\begin{equation}\label{eq:43}
Min\;z_{total}=\left[\frac{z_{1}-z_{1}^{*}}{z_{1}^{*}}+\frac{z_{2}-z_{2}^{*}}{z_{2}^{*}}+\frac{z_{3}^{*}-z_{3}}{z_{3}^{*}}\right]
\end{equation}

In equation (\ref{eq:43}), $z_1$ and $z_2$ are minimization objectives and their values will be greater than or equal to their optimal values, $z_1^*$ and $z_2^*$, respectively. On the other hand, $z_3$ is a maximization objective and its value will be smaller than or equal to its optimal value, $z_3^*$. As the deviation of an objective from its optimal value cannot be negative, equation (\ref{eq:43}) is written such that it remains positive $\left( z_1 - z_1^* \geq 0; z_2 - z_2^* \geq 0; z_3^* - z_3 \geq 0\right)$.

\section{Computational results}\label{Computational results}
In order to validate the described RO model for green supplier evaluation, we present results of various instances of the model. The experiments were conducted based on the data set of table \ref{tab: Datafor numexmple}, which was generated randomly using uniform distribution. Note that the parameters related to price were generated for the first period, and the next periods were calculated according to model constraints. Also, three different scenarios were considered, namely pessimistic, most realistic, and optimistic. 

\begin{table}[htbp]
\centering
\caption{Data for Numerical Example}
\scalebox{1}{
\renewcommand{\arraystretch}{0.5}
\begin{tabular}{|c|c|c|}

\hline
\textbf{Sets and} & \textbf{Scenarios} & \textbf{Amounts
} \\
\textbf{parameters} & & \\
\hline
$i$ & & $4$
 \\
\hline
$j$ & & $5$
 \\
\hline
$t$ & & $4$
 \\
\hline
$k$ & & $3$
 \\
\hline
$s$ & & $3$
 \\
\hline
$m$ & & $3$
 \\
\hline
$ir$ & & $0.04$
 \\
\hline
$\lambda _{1}$
 & & $15$
 \\
\hline
$\lambda _{2}$
 & & $15$
 \\
\hline
$\omega $
 & & $50$
 \\
\hline
$Pr^{s}$
 & $1$, 2, 3 & $0.2$, 0.6, 0.2 \\
\hline
$C_{ij1}^{s}$
 & $1$, 2, 3 & $U(10,23)$, $U(11.5,26)$, $U(13,30)$ \\
\hline
$h_{i1}$
 & & $U(28,35)$
 \\
\hline
$f_{i1}$
 & & $U(33,41)$
 \\
\hline
$L_{jt}^{s}$
 & $1,2,3$
 & $U(0,5)$
 \\
\hline
$G_{ijl}$
 & & $U(6,12)$
 \\
\hline
$e_{it}^{s}$
 & $1$, 2, 3 & $U(0.03,0.092)$, $U(0.035,0.126)$, $
U(0.04,0.145)$ \\
\hline
$p_{it}^{s}$
 & $1$, 2, 3 & $U(0.02,0.08)$, $U(0.023,0.092)$, $
U(0.027,0.105)$ \\
\hline
$u_{it}^{s}$
 & $1$, 2, 3 & $U(0.6,0.9), U(0.62,0.93)$, $U(0.63,0.94)$ \\
\hline
$v_{it}^{s}$
 & $1$, 2, 3 & $U(0.6,0.9), U(0.72,0.93), U(0.73,0.94)$
 \\
\hline
$O_{ij1}$
 & & $U(5,11)$
 \\
\hline
$SP_{t}^{m}$
 & & $U(4000,4020)$
 \\
\hline
$BP_{t}^{m}$
 & & U(3980,4000) \\
\hline
$DC_{i1}$
 & & U(4,7) \\
\hline
$RC_{i1}$
 & & U(10,17) \\
\hline
 $DP_{i1}$ & & U(3,5) \\
\hline
$TC_{j1kn}$
 & & $TC_{j111}=U(28,37),\; TC_{j112}=U(29,38) ,\;
TC_{j121}=U(35,40)$\\ & & $TC_{j122}=U(36,41)$,\; $TC_{j131}=U(39,52)$,\; $
TC_{j132}=U(40,53)$ \\
\hline
$DE_{it}^{s}$
 & $1, 2, 3$
 & $U(2500,4600)$, $U(2930,4760)$, $U(3070,4990)$ \\
\hline
$d_{j}$
 & & $U(3,7)$
 \\
\hline
$CET_jtkn$
 & & $CET_{jt1n}=U(0.29,0.37)$, $CET_{jt2n}=U(0.33,0.46)$ \\
& & $CET_{jt3n}=U(0.41,0.49)$
 \\
\hline
$CEP_{ijt}$
 & & $U(0.006,0.012)$
 \\
\hline
$CER_{ijt}$
 & & $U(0.006,0.012)$
 \\
\hline
$EM_{ijt}$
 & & $U(1,10)$
 \\
\hline
$GP_{ijt}$
 & & $U(1,10)$
 \\
\hline
$RE_{ijt}$
 & & $U(1,10)$
 \\
\hline
$PT_{ijt}$
 & & $U(1,10)$
 \\
\hline
$CAP_{t}$
 & & $U(170,200)$
 \\
\hline
$M_{k}$
 & & ${3000,6000,14000}$
 \\
\hline
\end{tabular}
\label{tab: Datafor numexmple}%
}
\end{table}

The results of the solved model are shown in table \ref{tab: Resultsnumexample}. The results determine the best suppliers and the optimal order allocations. For instance, $X_{112}=2617$ indicates that the manufacturer should place an order for 2617 of product 1 to the first supplier in the second period. Moreover, the optimal shortage and storage amounts, which are calculated based on the optimal orders are provided. The over-fulfillment and under-fulfillment values express the extra orders and the ignored demand, respectively. These amounts are calculated based on the relative importance the manufacturer assigns to solution optimality, i.e., achieving the best possible answer, compared to customer satisfaction. Furthermore, the optimum carbon trade quantity in each period is given. The results demonstrate that the optimal quantities for purchasing carbon allowance are 96.939, 64.684, 88.764, and 134.056 in periods 1, 2, 3, and 4, respectively. By using this data, managers will be able to find the right amount to order from each supplier considering the defined criteria, and how much carbon allowance they need to buy or sell in each period.

\begin{table}[h!]
  \centering
  \caption{Results of Numerical Example}
  \renewcommand{\arraystretch}{0.6}
    \begin{tabular}{|c|c|}
   \hline
\textbf{Variables} & \textbf{Results} \\
\hline
$X_{ijt}$
 & $X_{112}=2,617 \;\;\;X_{141}=1,825 \;\;\;X_{142}=2,506 \;\;\;X_{152}=1,627 \;$ \\ &
$X_{153}=4,091\;\;\;\; X_{211}=2,558\;\;\;\;
X_{212}=381 \;\;\;\;\;\;X_{214}=3,121$\; \\ & $X_{231}=126 \;\;\;\;\;\;\; X_{243}=810\;\;\; X_{252}=2,340 \;\;\;X_{253}=1,908$ \;\;\;\;\;\;\\ &
$X_{311}=2,962 \;\;\; X_{312}=2,756\;\;\; X_{343}=2,813  \;\;\;X_{412}=2,861$ \\ &
$X_{414}=2,690 \;\;\;X_{421}=3,382$
 \\
\hline
$r_{it}$
 & $0$
 \\
\hline
$b_{it}$
 & $b_{14}=4,243 \;\;\;b_{34}=3,056 \;\;\;b_{43}=3,657 \;\;\;b_{44}=3,937$
 \\
\hline
$\alpha _{t}$
 & $\alpha _{1}=96.939 \;\;\;\alpha _{2}=64.684\;\;\; \alpha _{3}=88.764\;\;\; \alpha 
_{4}=134.056$
 \\
\hline
$\beta _{t}$
 & $0$
 \\
\hline
$\delta _{ijt}^{s+}$
 & $\delta _{111}^{s+}=178\;\;\; \delta _{112}^{s+}=159 \;\;\;\delta _{124}^{s+}=202 
\;\;\;\delta _{153}^{s+}=147$ \\  &
$\delta _{214}^{s+}=123 \;\;\;\delta _{222}^{s+}=101 \;\;\;\delta _{243}^{s+}=106 
\;\;\;\delta _{251}^{s+}=92$ \\  &
$\delta _{312}^{s+}=97 \;\;\;\delta _{313}^{s+}=113 \;\;\;\delta _{341}^{s+}=126 \;\;\;\delta _{354}^{s+}=145$ \\  &
$\delta _{432}^{s+}=103 \;\;\;\delta _{441}^{s+}=143 \;\;\;\delta _{443}^{s+}=174 
\;\;\;\delta _{454}^{s+}=105$
 \\
\hline
$\delta _{ijt}^{s-}$
 & $\delta _{111}^{s-}=184\;\;\; \delta _{114}^{s-}=212 \;\;\;\delta _{123}^{s-}=150 
\;\;\;\delta _{142}^{s-}=164$ \\  &
$\delta _{221}^{s-}=94 \;\;\;\delta_{232}^{s-}=104 \;\;\;\delta 
_{253}^{s-}=109 \;\;\;\delta _{254}^{s-}=127$ \\  &
$\delta _{314}^{s-}=153 \;\;\;\delta _{321}^{s-}=131 \;\;\;\delta _{322}^{s-}=99 
\;\;\;\delta _{343}^{s-}=117$ \\  &
$\delta _{411}^{s-}=148 \;\;\;\delta _{412}^{s-}=105 \;\;\;\delta _{443}^{s-}=183 
\;\;\;\delta _{444}^{s-}=107$
 \\
\hline
$z_{1}$
 & $640,743$
 \\
\hline
$z_{2}$
 & $43,627$
 \\
\hline
$z_{3}$
 & $1,133,472$
 \\
\hline
$z_{Total}$
 & $1.699$
 \\
\hline
    \end{tabular}%
  \label{tab: Resultsnumexample}%
\end{table}%
\vspace{1cm}

To analyze the model and check the sensitivity of its objectives to parameters, several analyses were performed. The results will be discussed in the following sections.

\subsection{Sensitivity analysis on $\omega$} \label{Sensitivity analysis on omega}
One of the important parameters in robust models is $\omega$, the value that weighs the trade-off between optimality and feasibility. In this model, infeasibility is defined as over-fulfillment or under-fulfillment. Over-fulfillment elevates the warehousing costs, while under-fulfillment results in lower customer satisfaction; customer's order is affected by the possibility of under-fulfillment. To decrease the risk of infeasibility, one should increase the penalty of receiving an infeasible solution, which, in turn, leads to reduced model optimality. Therefore, there is always a trade-off between model infeasibility and optimality. Understandably, it is not easy for decision-makers to determine the exact value of $\omega$. Consequently, a suitable way to investigate this issue is to measure different values of the objectives and the amount of model infeasibility when distinct values for $\omega$ are considered. Afterward, an analysis of the results enables the decision-makers to manipulate the preferred value of $\omega$.

\begin{figure}[h]
\includegraphics[width=\textwidth]{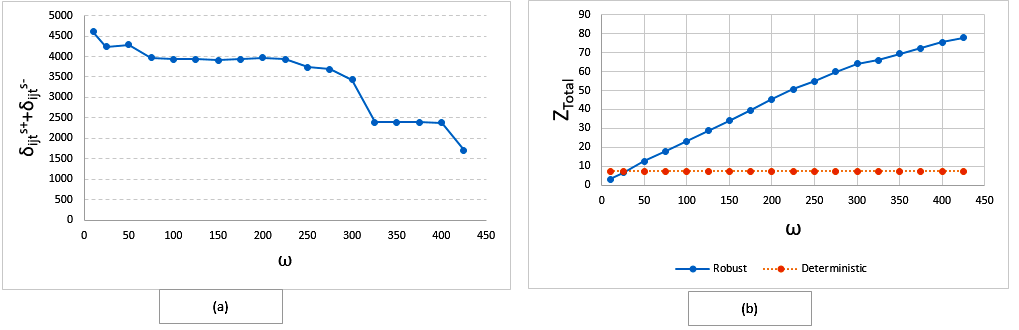}
\caption{Sensitivity analysis on $\omega$}
\label{Fig: Sensitivity analysis on omega}
\end{figure}

By increasing the value of $\omega$, due to the imposed penalty, the model tries to decrease the infeasibility. Therefore, as shown in figure \ref{Fig: Sensitivity analysis on omega}(a), when $\omega$ increases, infeasibility, which is defined as the amount of over-fulfillment and under-fulfillment, decreases. On the other hand, by increasing the value of $\omega$, firms will be penalized more for infeasibility; in this situation, first and second objective functions, which are formulated as robust equations, will deteriorate, and consequently $z_{Total}$ will become worse. This is shown in figure \ref{Fig: Sensitivity analysis on omega}(b).

By scrutinizing these results, decision-makers can assign a precise value to $\omega$ based on the goals defined by the firm. For instance, if the objective values are more important, a lower weight for model robustness should be set. On the other hand, if customer satisfaction is more critical, more penalty must be in place for infeasibility (higher values for $\omega$), which leads to worse objective value.

Furthermore, in figure \ref{Fig: Sensitivity analysis on omega}(b), a comparison is performed between the values of robust and deterministic models. The figure depicts that for a large domain of $\omega$ values ($\omega > 25$) the objective function of the deterministic model is better than the robust model. This occurs due to the penalty imposed to the robust model for deviation and infeasibility, i.e., there is more restriction on the objective function. Regardless, by handling uncertainty, robust models are able to achieve solutions that are closer to real world scenarios. This enhances their applicability.

\subsection{Sensitivity analysis on $\lambda$}\label{Sensitivity analysis on lambda1 and lambda2}
\begin{figure}[h]
\includegraphics[width=\textwidth]{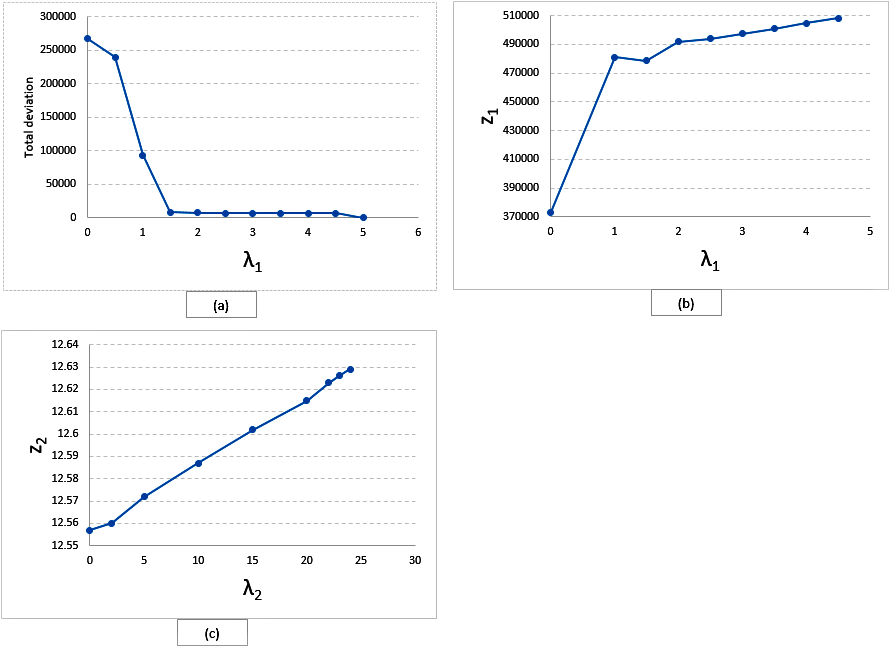}
\caption{Sensitivity analysis on $\lambda$}
\label{Fig: Sensitivity analysis on lambda}
\end{figure}

Another crucial aspect is the analysis of the impact of the penalty on the deviation of scenarios from the expected value of the objective. Note that the model considers various scenarios, and one cannot be certain which one will happen. Hence, the solution method seeks to achieve a solution which is the closest possible value to the optimum of all scenarios. In this sense, $\lambda$ is the weight parameter that alludes to the importance of solution closeness to the real-world application. A larger value for $\lambda$ implies a higher level of importance for deviation from the average; plus, the objective value deteriorates due to the increased penalty. This is illustrated in figure \ref{Fig: Sensitivity analysis on lambda}. As displayed in figure \ref{Fig: Sensitivity analysis on lambda}(a), by increasing the value of $\lambda_{1}$, the deviation between the objective value in each scenario and the average objective value decreases. On the other hand, according to figure \ref{Fig: Sensitivity analysis on lambda}(b), the value of the first objective function increases. In other words, the first objective function value exacerbates as it is a minimization objective. Consequently, total objective value deteriorates. 

The same occurs for $\lambda_{2}$, but since in this paper the values of the second objective are significantly smaller than the first objective, $\lambda_{2}$ has negligible effect on the deviation of the objective function (see table \ref{tab: landa table}). It does not mean that the second objective is not important, because the total objective is normalized using equation \eqref{eq:43} to avoid overestimating the importance of objectives with large values. However, the deviation is fixed on a small amount (11.148) compared to the first objective (486,736), and adding the value of $\lambda_{2}$ does not have an impact on the deviation, although the second objective function worsens (figure \ref{Fig: Sensitivity analysis on lambda}(c)). Therefore, it is recommended to set the value of $\lambda_{2}$ to a small number to avoid deteriorating the second objective, which represents the amount of emission.

In other words, putting more penalty on the deviation of objective values in different scenarios reduces the deviation, but worsens the objective values. This provides a road map for the decision-makers to choose values of $\lambda_{1}$ and $\lambda_{2}$ according to their goals.

\begin{table}
\begin{center}

\centering
\caption{Effect of $\lambda_{2}$ on Deviation of Second Objective Value}
\begin{tabular}{|c|c|c|c|c|c|c|} 
\hline
 {$\lambda_{2}$} & $z_{Total}$ & $z_{1}$      & $z_{2}$     & $z_{3}$        & $\delta$     & Deviation from Second Objective  \\ 
\hline
{0}                                    & {12.557} & {486,736} & {43,182} & {1,170,111} & {4,283} & {11.148}                                                                                                 \\
\hline
{2}                                    & {12.56}  & {486,736} & {43,193} & {1,170,111} & {4,283} & {11.148}                                                                                                 \\
\hline
{5}                                    & {12.572} & {486,736} & {43,237} & {1,170,111} & {4,283} & {11.148}                                                                                                 \\\hline
{10}                                   & {12.587} & {486,736} & {43,291} & {1,170,111} & {4,283} & {11.148}                                                                                                 \\\hline
{15}                                   & {12.602} & {486,736} & {43,346} & {1,170,111} & {4,283} & {11.148}                                                                                                 \\\hline
{22}                                   & {12.623} & {486,736} & {43,422} & {1,170,111} & {4,283} & {11.148}                                                                                                 \\\hline
{23}                                   & {12.626} & {486,736} & {43,433} & {1,170,111} & {4,283} & {11.148}                                                                                                 \\\hline
{24}                                   & {12.629} & {486,736} & {43,444} & {1,170,111} & {4,283} & {11.148}                                                                                                 \\
\hline
\end{tabular}
\label{tab: landa table}%
\end{center}
\end{table}

\subsection{Sensitivity analysis on CAP}\label{Sensitivity analysis on CAP}
As mentioned in section \ref{Introduction}, governments usually reduce the assigned cap in cap-and-trade mechanism each year to decrease air pollution\footnote{As an example, see the cap-and-trade mechanism of Ontario, Canada (accessed March 13, 2020): \url{https://www.ontario.ca/page/cap-and-trade}.}. Therefore, an analysis of cap reduction is required to predict the probable effects of this restriction and make the right actions toward the change. 

\begin{figure}[h]
\includegraphics[width=\textwidth]{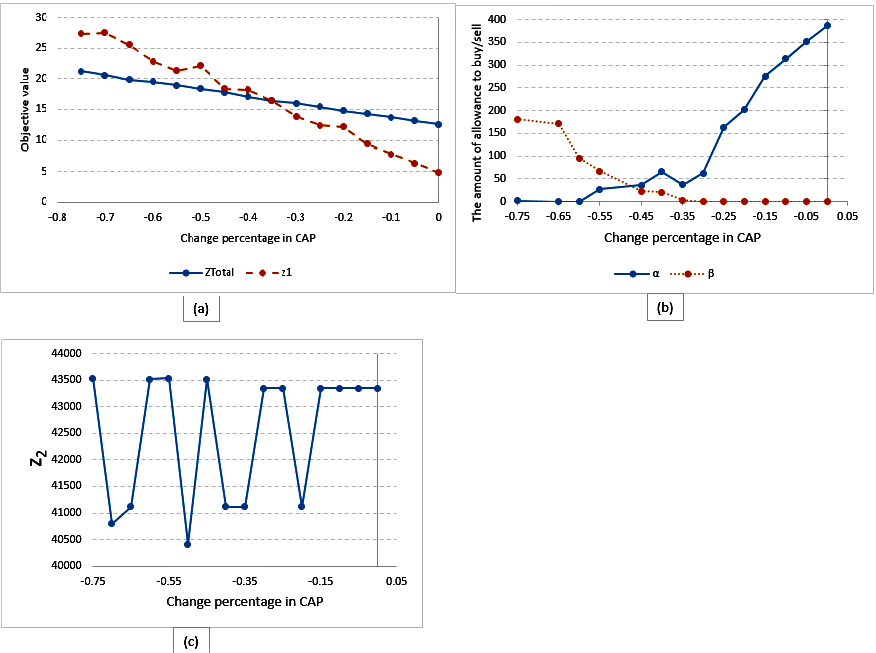}
\caption{Sensitivity analysis on CAP$^\dagger$\\$^\dagger$Values on the x-axis show cap reduction. For example, -0.1 means cap is decreased by 10\%.}
\label{Fig: Sensitivity analysis on CAP}
\end{figure}

Figure \ref{Fig: Sensitivity analysis on CAP}(a) shows a reverse relationship between objective values and carbon cap. It indicates that by decreasing the carbon allowance cap, the cost objective function ($\times 10^{5}$) increases, thus worsening the total objective value. This matter can be investigated from different perspectives. From government's point of view, the total utility decreases by decreasing the cap, but it results in less carbon emission, which can be deemed as a positive result for government due to health- and environmental-related costs. From the manufacturer's point of view, decreasing the cap is not an ideal option, because they are forced to purchase more allowance, or they will make less money for not being able to sell as much allowance, or they will have to upgrade to more expensive green technologies.

In addition, this restriction affects the amount of allowance to buy or sell. As displayed in figure \ref{Fig: Sensitivity analysis on CAP}(b), while the cap decreases, the amount of total carbon allowance the manufacturer needs to buy (sell), increases (decreases). This is the main driver for the increased cost. 

Furthermore, according to figure \ref{Fig: Sensitivity analysis on CAP}(c), changing the cap does not have a discernible impact on the optimal amount of carbon emission ($z_{2}$); no trend can be verified in this figure. For the parameter values featured in this study, a cap reduction does not significantly affect the carbon released by manufacturers since they consider the carbon allowance market as a source of profit. In other words, manufacturers focus on the cost function while trying to satisfy the restrictions related to carbon emission. 

Using these results help the government assign the proper carbon emission cap and minimize the impact on manufacturers in terms of emission costs.

\subsection{Sensitivity analysis on cap-and-trade market prices}\label{Sensitivity analysis on cap and trade market prices}
Buying and selling prices are critical parameters in the cap-and-trade system with a direct impact on the amount of carbon allowance being traded by firms. According to table \ref{tab: SensitivityonBPSP}, by increasing the selling price, the values of the first and total objectives decrease. This can be traced to the growth in the optimal value of the amount of sold allowance.

As shown in table \ref{tab: SensitivityonBPSP}, when $BP > SP$ the manufacturers turn into brokerages. In this case, firms can buy allowance in the market and sell it at a higher price for a profit. Therefore, the optimal values of variables that control buying and selling allowance in the market will be positive in the same period. To avoid this situation, government can restrict the carbon market price such that $BP$ is always less than $SP$ in carbon trading market.

\begin{table}
\centering
\caption{Results of sensitivity analysis on carbon market prices}
\begin{tabular}{|c|c|c|c|c|c|c|c|} 
\hline
 Parameter & Price change (\%) & {$z_{Total}$} & {$z_{1}$}      & {$z_{2}$}     & {$z_{3}$}        & {$\alpha$}       & {$\beta$}        \\ 
\hline
{}{}{{BP}}                & {-20}      & {13.24}  & {847,973} & {40,670} & {1,365,767} & {296}     & {0}        \\
                                            & {-10}      & {12.783} & {728,194} & {40,670} & {1,365,767} & {296}     & {0}        \\
                                            & {-5}       & {12.9}   & {555,835} & {43,466} & {1,147,220} & {377}     & {0}        \\
                                            & {0}        & {12.602} & {486,736} & {43,346} & {1,170,111} & {387}     & {0}        \\
                                            & {5}        & {10.108} & {39,832}  & {40,497} & {1,396,833} & {5,343}   & {5,085}    \\
                                            & {10}       & {-0.306} & {1,363}   & {3,684}  & {5,571,751} & {3.9×108} & {3.9×108}  \\ 
\hline
{}{}{\textit{SP}}                & {-10}      & {-0.307} & {1,363}   & {3,681}  & {5,568,767} & {9.9×107} & {9.9×107}  \\
                                            & {-5}       & {-0.307} & {1,363}   & {3,688}  & {5,605,797} & {135,493} & {136,376}  \\
                                            & {0}        & {12.602} & {486,736} & {43,346} & {1,170,111} & {387}     & {0}        \\
                                            & {5}        & {12.602} & {486,736} & {43,346} & {1,170,111} & {385}     & {0}        \\
                                            & {10}       & {12.602} & {486,736} & {43,346} & {1,170,111} & {387}     & {0}        \\
                                            & {20}       & {12.602} & {486,736} & {43,346} & {1,170,111} & {387}     & {0}        \\
\hline
\end{tabular}
\label{tab: SensitivityonBPSP}%
\end{table}

\subsection{Cap-and-trade versus paying penalty to the government}\label{Cap and trade vs. paying the penalty to the government}
\begin{figure}[h]
\includegraphics[width=\textwidth]{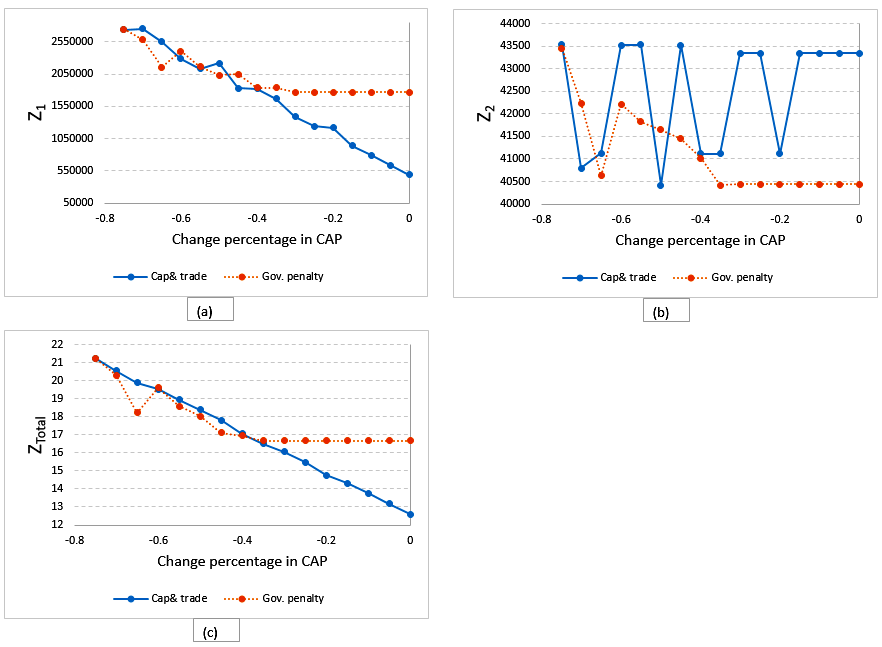}
\caption{Cap-and-trade vs. penalty-based}
\label{Fig: Cap and trade vs. penalty-based}
\end{figure}

Another regulation possibility to investigate is forcing the manufacturers to pay a penalty to the government for carbon emission once a certain cap is exceeded. This policy does not allow formation of a carbon emission trading market. As shown in figures \ref{Fig: Cap and trade vs. penalty-based}(a)-\ref{Fig: Cap and trade vs. penalty-based}(c), for the cap changes below 40 percent, cost of pollution is less in cap-and-trade compared to the penalty-based mechanism. However, the amount of carbon released under cap-and-trade policy is more than the emissions in the government penalty regulation. Therefore, in case there is a high priority on the carbon reduction objective, government must decrease the cap in the C\&T mechanism to a value which is less than 40\% of the original amount. Also, $z_{Total}$ in cap-and-trade is lower in the given range. In other words, cap-and-trade performs better than the penalty-based system for manufacturers. Namely, cap-and-trade derives a better total utility; thus, it is preferred to the penalty-based approach.

If the cap is lowered more than 40 percent, no superiority exists between the two mechanisms. This occurs since the cap is decreased below the 40 percent mark, the amount of allowance the manufacturer needs to buy (or the amount of penalty the manufacturer should pay to the government) will be positive; and the quota to sell will gradually tend to zero. In other words, the manufacturer will not be able to generate profit by selling emission allowance. Hence, the gap between the cap-and-trade mechanism and the penalty system will vanish. The above scrutiny is an important result for the governments to assign effective values for the cap according to the defined emission targets, while not hurting the businesses. 

Based on the presented results, cap-and-trade culminates in better outcomes in terms of total model utility, although with higher carbon emission which is the result of assigning equal weights to economic selection criteria, green selection criteria, and carbon emission level. Obviously, assigning more weight to carbon emission level translates to solutions with less carbon emission.

Another inquiry on the two carbon control systems is based on the fact that cap-and-trade mechanism allocates a constant carbon emission permission to the entire market. Therefore, the manufacturers buy the allowance of the other players in the carbon trading market, while the total allowed carbon emission is a fixed value. AS a result, if the firms choose not to sell their allowance in the market, or if they consume all of their own allowance, other firms will not be able to place orders anymore. 

Conversely, in penalty-based systems, although there is a cap determined by the government, the total amount of carbon released can exceed the cap because manufacturers can emit as much carbon as they need by paying the penalty to the government. Therefore, the total carbon released by manufacturers may be more in this system. Note that this argument cannot be confirmed by the results of this paper because data from all of the manufacturers in the market is required to audit the total amount of carbon emission. Consequently, this topic is left for study in future research efforts.

\section{Conclusion}\label{Conclusion}
To support decisions related to outsourcing in a green supply chain and to help with procuring the best materials from the best suppliers in an uncertain environment encompassing both green and traditional criteria, a multi-objective robust optimization model is generalized in this paper. The cap-and-trade mechanism is considered to evaluate the role of this pollution control system. A numerical example is presented to analyze different aspects of the generalized model and the solution approach. 

The optimal activity levels of a firm to achieve the maximum profit while meeting the environmental goals are obtained. The proposed model enables the decision makers to work with a more flexible decision support system due to the analysis done on $\omega$ and $\lambda$. The results show the effect of cap amount and the trade prices on the firm's objective. Therefore, manufacturers will be well-informed about selecting the best suppliers, the amount of orders to place with each of the selected suppliers, and trading in the allowance market to meet their objectives. Additionally, the superiority of cap-and-trade mechanism compared to penalty-based system regarding the total utility of supply chain from a micro-economic perspective is shown in this paper.

Since this paper is one of the first to consider cap-and-trade mechanisms in the green supply chain, more research in this area can be pursed. For instance, this paper analyzes the cap changes and its effects on the objectives from the manufacturer's point of view. Accordingly, a motivating topic for future research is finding the optimal value of the cap assigned by the government to reduce the adverse environmental effects of production industries, while remaining business friendly. One could explore a pricing study to find a reasonable range for the allowance prices according to different factors. Specifically, decreasing the cap will leave all of the manufacturers with less allowance for carbon emission. Therefore, the manufacturers will have less allowance to offer in the trade market, resulting in elevated trading prices. In this case, some manufacturers may decide to simply sell their allowance to generate profit. This presents a challenge from the existing problem’s perspective.

In this paper, we studied the difference between cap-and-trade system and penalty-based mechanism to control the air pollution from the manufacturer's point of view. If one has access to data about the total amount of carbon released by all manufacturers in the market, one could study the differences between the two particular pollution control regimes from the government's perspective.

\bibliography{References}

\end{document}